\documentclass[11pt]{article}
\usepackage[T1]{fontenc}
\usepackage{lmodern}
\usepackage[a4paper,margin=25mm]{geometry}
\usepackage{amsmath,amssymb,amsthm,mathtools}
\usepackage{graphicx,pgf,booktabs,longtable,array}
\usepackage{pdflscape}
\usepackage{microtype}
\usepackage[colorlinks=true,linkcolor=blue!45!black,citecolor=blue!45!black,urlcolor=blue!45!black]{hyperref}
\usepackage{enumitem}
\usepackage{placeins}

\newif\ifrevisionnotes
\revisionnotestrue 
\newtheorem{theorem}{Theorem}[section]
\newtheorem{proposition}[theorem]{Proposition}
\newtheorem{lemma}[theorem]{Lemma}

\theoremstyle{remark}

\numberwithin{equation}{section}
\DeclareMathOperator{\Li}{Li}

\DeclareMathOperator{\arccosh}{arccosh}

\DeclareMathOperator{\Var}{Var}

\DeclareMathOperator{\arsinh}{arsinh}

\newcommand{\C}{\mathbb C}
\newcommand{\R}{\mathbb R}
\newcommand{\E}{\mathbb E}
\newcommand{\Prob}{\mathbb P}
\newcommand{\dd}{\,\mathrm d}
\newcommand{\dA}{\,\mathrm dA}

\newcommand{\ol}[1]{\overline{#1}}
\title{Random complex cross ratios}
\author{Tim Marshall \and Gaven Martin \thanks{The research of both authors is supported by the New Zealand Marsden Fund.\newline Institute for Advanced Study, Massey University.}}
\begin{document}
\maketitle
\begin{abstract}
We study the cross ratio of four independent points chosen uniformly on the
Riemann sphere. A positive integral over a two-dimensional simplex gives its
normalized planar density. We then
evaluate this integral in terms of the Bloch--Wigner dilogarithm, providing a
closed formula. We
describe its symmetries, real-axis values, moment bounds and geometric
interpretations. We use these to evaluate the mean
ideal-tetrahedron volume directly from the density, recovering $\pi/6$,
a recent theorem of Kabluchko and Schange. We then prove
that the expected distance between two independently sampled hyperbolic geodesics in $\mathbb{H}^3$ is
$4/3$, by evaluating the logarithmic potential of the cross-ratio
distribution on the negative real axis. Further applications give the joint
dihedral-angle density and an exact acute-tetrahedron probability, and a
branched-cover probability law on the modular surface with exponentially
small cusp tails.
\end{abstract}

\section{The probability model}
Let $\widehat\C$ denote the Riemann sphere and $\mathrm dA(z)=\mathrm dx\,\mathrm dy$
the spherical area. Normalized
spherical area is the probability measure
\begin{equation}\label{eq:spherical}
 \dd\mu(z)=p(z)\dA(z),\qquad
 p(z)=\frac{1}{\pi(1+|z|^2)^2}.
\end{equation}
Throughout, $Z_1,Z_2,Z_3,Z_4$ are independent with law $\mu$, and our cross-ratio
convention is
\begin{equation}\label{eq:crossratio}
 Z=[Z_1,Z_2,Z_3,Z_4]
   =\frac{(Z_1-Z_3)(Z_2-Z_4)}{(Z_1-Z_2)(Z_3-Z_4)}.
\end{equation}
Coincident points have probability zero. All densities below are with respect
to $\mathrm dA$, unless we explicitly identify another measure.

\medskip

The moduli space $M_{0,4}$ and $\C\setminus\{0,1\}$ identify via
cross ratio. Forgetting the labels takes the quotient by the six anharmonic
symmetries. The measure studied here is then the pushforward of four copies
of normalized spherical area. This extends the circle model considered earlier in \cite{Martin}.

\medskip
Section~\ref{sec:applications} develops three consequences of the density.
The double cover branched over the four points provides a distribution on
elliptic curves and punctured tori, which we compare with the normalized
hyperbolic-area model of \cite{MartinLattices}. This construction is distinct
from the side-pairing construction of \cite{Martin}. A further application
to commutator parameters is motivated by the random-group programme of
\cite{MO} and so its precise sampling model is specified below.

For a sample $W\sim\mu$, its argument is uniform and independent of its modulus
$R$, whose density and distribution function are
\begin{equation}\label{eq:radial}
 p_R(r)=\frac{2r}{(1+r^2)^2},\qquad
 \Prob(R\leq r)=\frac{r^2}{1+r^2},\qquad r\geq0.
\end{equation}
Thus, for independent $U\sim\mathrm{Unif}(0,1)$ and
$\Theta\sim\mathrm{Unif}(0,2\pi)$,
\begin{equation}\label{eq:sampler}
 W=\sqrt{U^{-1}-1}\,e^{i\Theta}
\end{equation}
has exactly the required law. The spherical isometry $w\mapsto1/w$ preserves
this law, so $W\overset{d}=1/W$. These two random variables are, of course,
dependent when formed from the same sample.

\section{A positive formula for the density}
A rotation taking $Z_1$ to $\infty$ preserves the conditional joint
law of the other three points: conditional on $Z_1$, they remain independent
with law $\mu$. Since cross ratio is invariant under that rotation,
\begin{equation}\label{eq:threepoints}
 Z\overset{d}=\frac{A-C}{B-C},\qquad A,B,C\text{ independent with law }\mu.
\end{equation}
Only spherical rotations are used to normalize the sampling law as typically 
M\"obius maps preserve cross ratio but do not preserve $\mu$.

Let
\[
 \Delta=\{(\alpha,\beta):\alpha>0,\ \beta>0,\ \alpha+\beta<1\},
 \qquad \gamma=1-\alpha-\beta.
\]
\begin{theorem}[Normalized simplex formula]\label{thm:simplex}
For $z\in\C\setminus\{0,1\}$, the cross-ratio density is
\begin{equation}\label{eq:simplex}
 \boxed{\displaystyle
 \rho(z)=\frac{2}{\pi}\int_\Delta
 \frac{\alpha\beta\gamma}
 {\bigl(\alpha\beta|1-z|^2+\alpha\gamma|z|^2+\beta\gamma\bigr)^2}
 \dd\alpha\dd\beta.}
\end{equation}
It is positive and real analytic away from $0$ and $1$, and
$\int_\C\rho(z)\dA(z)=1$.
\end{theorem}
\begin{proof}
Put $C=c$, $B=c+v$ and $A=c+zv$. The Jacobian of the map
$(z,v,c)\mapsto(A,B,C)$ is $|v|^2$. Hence
\begin{equation}\label{eq:complexintegral}
 \rho(z)=\frac{1}{\pi^3}\int_{\C^2}
 \frac{|v|^2\dA(c)\dA(v)}
 {(1+|c+zv|^2)^2(1+|c+v|^2)^2(1+|c|^2)^2}.
\end{equation}
The elementary parameter identity
\[
 \frac{1}{P^2Q^2R^2}
 =120\int_\Delta\frac{\alpha\beta\gamma}
 {(\alpha P+\beta Q+\gamma R)^6}\dd\alpha\dd\beta
 \qquad(P,Q,R>0)
\]
follows by writing each reciprocal square as a gamma integral and separating
the sum of its three parameters from their proportions. All integrands are
nonnegative, so Tonelli's theorem justifies interchanging integrals.

Set $q=\alpha\beta|1-z|^2+\alpha\gamma|z|^2+\beta\gamma$.
Completing the square gives
\[
 1+\alpha|c+zv|^2+\beta|c+v|^2+\gamma|c|^2
 =1+|c+(\alpha z+\beta)v|^2+q|v|^2.
\]
The two complex integrations are now radial:
\[
 \int_{\C^2}\frac{|v|^2\dA(c)\dA(v)}
 {(1+|c|^2+q|v|^2)^6}
 =\frac{\pi}{5}\int_\C\frac{|v|^2\dA(v)}{(1+q|v|^2)^5}
 =\frac{\pi^2}{60q^2}.
\]
Substitution proves \eqref{eq:simplex}, with normalization inherited from
the change of variables in \eqref{eq:complexintegral}.
On a compact subset of $\C\setminus\{0,1\}$, $q$ is bounded below by a
positive constant times $\alpha\beta+\alpha\gamma+\beta\gamma$.
The resulting dominating integrand is integrable at the three vertices of
$\Delta$. A locally uniform power-series expansion of $q^{-2}$ therefore
also proves real analyticity.
\end{proof}

\subsection{A stable one-dimensional quadrature}
Set $a=|z|^2$, $b=|1-z|^2$, and, for $0<u<1$, put
\[
 U=u(1-u),\qquad L=1+(a-1)u,\qquad M=bU.
\]
Using $\alpha=su$, $\beta=s(1-u)$, $\gamma=1-s$ in
\eqref{eq:simplex} and integrating $s$ gives
\begin{equation}\label{eq:oned}
 \rho(z)=\frac{2}{\pi}\int_0^1 U\,K(L,M)\dd u,
 \qquad
 K(L,M)=\frac{(L+M)\log(L/M)-2(L-M)}{(L-M)^3}.
\end{equation}
Positivity is immediate from the equivalent expression
\begin{equation}\label{eq:kernel}
 K(L,M)=\int_0^1\frac{s(1-s)}{((1-s)L+sM)^2}\dd s.
\end{equation}
The apparent singularity at zero when $L=M$ is removable, with $K(M,M)=1/(6M^2)$.
For numerical evaluation near  $L=M$ we found it best to use 
\begin{equation}\label{eq:series}
 K(M(1+\delta),M)=\frac{1}{M^2}
 \sum_{j=0}^{\infty}\frac{(-1)^j(j+1)}{(j+2)(j+3)}\delta^j,
 \qquad |\delta|<1.
\end{equation}
as this gives an check of the closed form and avoids
cancellation near the real axis.

\section{Evaluation by the Bloch--Wigner dilogarithm}
The Bloch--Wigner function is
\begin{equation}\label{eq:BW}
 D(z)=\operatorname{Im}\Li_2(z)+\arg(1-z)\log|z|,
 \qquad \Li_2(z)=-\int_0^z\frac{\log(1-t)}{t}\dd t,
\end{equation}
with the principal branches off their cuts and the continuous extension
$D(x)=0$ for real $x$. It satisfies \cite{Zagier}
\begin{align}\label{eq:Dsym}
 D(\ol z)&=-D(z),\nonumber\\
 D(z)&=-D(1-z)=-D(1/z)=-D(z/(z-1))\nonumber\\
     &=D(1-1/z)=D(1/(1-z)).
\end{align}
We will use the differential identity
\begin{equation}\label{eq:Ddiff}
 \mathrm dD(z)=\log|z|\,\mathrm d\arg(1-z)
              -\log|1-z|\,\mathrm d\arg z,
\end{equation}
which follows directly by differentiating \eqref{eq:BW}.

\subsection{The master integral}
For independent positive parameters $a,b$, define
\begin{equation}\label{eq:J}
 J(a,b)=\int_\Delta\frac{\dd\alpha\dd\beta}
 {a\alpha\gamma+b\alpha\beta+\beta\gamma}
 =\int_0^1\frac{\log(L/M)}{L-M}\dd u.
\end{equation}
The second equality uses the substitution above.

\begin{lemma}\label{lem:master}
If $z=x+iy$ with $y>0$, $a=|z|^2$ and $b=|1-z|^2$, then
\begin{equation}\label{eq:master}
 J(a,b)=\frac{2D(z)}{y}.
\end{equation}
\end{lemma}
\begin{proof}
Substitute $u=t/(1+t)$ into \eqref{eq:J}. Since $1+a-b=2x$,
\[
 J(a,b)=\int_0^\infty
 \frac{\log\bigl((1+t)(1+at)/(bt)\bigr)}{|1+zt|^2}\dd t.
\]
Consider the two paths
\[
 w_1(t)=\frac{z(1+t)}{1+zt},\qquad
 w_2(t)=\frac{1+zt}{1+at}.
\]
Both lie in the upper half-plane for $0<t<\infty$.
Equation \eqref{eq:Ddiff} now gives
\begin{align*}
 \frac{\mathrm d}{\mathrm dt}D(w_1(t))
 &=\frac{y}{|1+zt|^2}
   \left(\frac12\log\frac{b}{a}-\log(1+t)\right),\\
 \frac{\mathrm d}{\mathrm dt}D(w_2(t))
 &=\frac{y}{|1+zt|^2}
   \left(-\log t-\frac12\log(ab)+\log(1+at)\right).
\end{align*}
Subtracting gives $y$ times the integrand of $J$.
At $t=0$, $(w_1,w_2)=(z,1)$, while at infinity they tend to
$(1,1/\ol z)$. By \eqref{eq:Dsym}, $D(1/\ol z)=D(z)$.
Integration therefore yields $yJ=2D(z)$.
\end{proof}

\subsection{Differentiation and the closed formula}
We now hold $a$ fixed. For
$g=\log(L/M)/(L-M)$ differentiation using $M=bU$ gives
\begin{equation}\label{eq:operator}
 \partial_b g+b\partial_b^2g
 =U\left[\frac{(L+M)\log(L/M)}{(L-M)^3}
          -\frac{2}{(L-M)^2}\right]=U K(L,M).
\end{equation}
This identity extends continuously through $L=M$. Differentiation under
\eqref{eq:J} is easily justified. Consequently,
\begin{equation}\label{eq:rhooperator}
 \rho(z)=\frac{2}{\pi}\,\partial_b\bigl(b\partial_bJ(a,b)\bigr).
\end{equation}
This is the useful link between a positive probability integral and its
dilogarithmic evaluation.

\begin{theorem}[Closed density formula]\label{thm:closed}
For $z=x+iy$ with $y\ne0$,
\begin{align}\label{eq:closed}
 \rho(z)={}&
 \frac{3x^2-3x^3+y^2-3xy^2}{2\pi y^4}\log|z|^2 +\frac{-6x^2+3x^3-2y^2+3x(1+y^2)}{2\pi y^4}\log|1-z|^2\nonumber\\
 &-\frac{1}{\pi y^2}
 +\frac{3(x-1)^2x^2+(2x-1)^2y^2+y^4}{\pi y^5}D(z).
\end{align}
\end{theorem}
\begin{proof}
For $y>0$, write $x=(a+1-b)/2$ and $y=(a-x^2)^{1/2}$.
Then $\partial_bx=-1/2$ and $\partial_by=x/(2y)$. From
\eqref{eq:Ddiff},
\[
 \partial_bD(z)=\frac{a-x}{4yb}\log a-\frac{1}{4y}\log b,
\]
so Lemma~\ref{lem:master} gives
\[
 \partial_bJ
 =\frac{(a-x)\log a/b-\log b}{2y^2}-\frac{xD(z)}{y^3}.
\]
One further differentiation in \eqref{eq:rhooperator} gives
\eqref{eq:closed}. Complex conjugation establishes the lower half-plane
case. 
\end{proof}

\section{Symmetries, special values and tails}
Permutation of the four independent points and conjugation imply that
\begin{equation}\label{eq:symmetry}
 \boxed{\rho(\ol z)=\rho(z)=\rho(1-z)
       =|z|^{-4}\rho(1/z).}
\end{equation}
More generally $\rho(z)=\rho(T(z))|T'(z)|^2$ for each of
\[
 T(z)\in\left\{z,1-z,\frac1z,\frac{z}{z-1},
                  1-\frac1z,\frac1{1-z}\right\}.
\]
The Jacobian here is necessary as invariance of the probability measure does imply
mean invariance of the planar density under inversion.

\subsection{Values on the real line}
The apparent poles in powers of $y$ in \eqref{eq:closed} cancel. For
$x\in\R\setminus\{0,1\}$, its continuous extension is
\begin{align}\label{eq:real}
 \rho(x)=\frac{2}{15\pi x^3(x-1)^3}\bigl[&-2x(x-1)(x^2-x+1)\nonumber\\
 &+x^3(2x^2-5x+5)\log|x|\nonumber\\
 &-(x-1)^3(2x^2+x+2)\log|x-1|\bigr].
\end{align}
For example,
\[
 \rho(-1)=\rho(2)=\frac{2\log2-1}{5\pi}
       =0.024592\ldots,\qquad
 \rho(1/2)=16\rho(-1).
\]
Formula \eqref{eq:real} follows either by expanding \eqref{eq:closed} at
$y=0$ or by integration in \eqref{eq:oned}. Note that this is the
\emph{trace of a two-dimensional density} and not the similar density for four points
on a circle determined in \cite{Martin}.

\subsection{Two checks on special curves}
On $z=e^{i\theta}$, with $x=\cos\theta$, $y=\sin\theta\ne0$,
\begin{equation}\label{eq:circle}
 \rho(z)=\frac{2(1-x+x^2)D(z)
 -y\{1+x+(1-2x)\log[2(1-x)]\}}{\pi y^3(1+x)}.
\end{equation}
On $z=1/2+iy$, $y\ne0$,
\begin{equation}\label{eq:midline}
 \rho(z)=\frac{-16y^3+(3+16y^4)D(z)
              +4y(3-4y^2)\log|z|}{16\pi y^5}.
\end{equation}
Both are direct substitutions in \eqref{eq:closed}; the coefficient of
$D(z)$ in \eqref{eq:midline} has a plus sign.

\begin{proposition}[Moment range]\label{prop:moments}
For real $p$,
\begin{equation}\label{eq:moments}
 \E|Z|^p<\infty\quad\Longleftrightarrow\quad -2<p<2.
\end{equation}
In particular, $\E Z=1/2$, while $\E|Z|^2=\infty$.
\end{proposition}
\begin{proof}
Here $f\asymp g$ means the ratio is bounded above and below by
positive constants. So, 
\[
 \rho(z)\asymp 1+\log(1/|z|)\quad(z\to0).
\]
In \eqref{eq:simplex}, let $\alpha=(1-\beta)v$ and
$\gamma=(1-\beta)(1-v)$. The integrand becomes
\[
 \frac{\beta(1-\beta)v(1-v)}
 {\{a(1-\beta)v(1-v)+\beta[bv+1-v]\}^2},
 \quad 0<\beta,v<1.
\]
As $z\to0$, $a=|z|^2\to0$ and $b$ stays bounded above and below away
from zero. On $\beta\leq1/2$, integrating the comparable expression
$\beta v(1-v)/(a v(1-v)+\beta)^2$ in $\beta$ gives an upper bound
of order $v(1-v)\big(1-\log[a v(1-v)]\big)$. This is integrable in $v$,
and restricting $v$ to $[1/4,3/4]$ gives a matching logarithmic lower bound.
The part where $\beta\geq1/2$ is uniformly bounded.

The symmetries now give a logarithmic singularity at $1$ and
$\rho(z)\asymp |z|^{-4}\log|z|$ at infinity. Polar integration proves
\eqref{eq:moments}, including divergence at both endpoints. Since the
first absolute moment is finite, $Z\overset d=1-Z$ implies
$\E Z=1-\E Z=1/2$.
\end{proof}

\subsection{The induced conformal metric and its punctures}\label{sec:metric}
The invariance of cross ratio under  conformal automorphisms
of the sphere identifies its parameter space with
$M_{0,4}\cong X:=\widehat\C\setminus\{0,1,\infty\}$.
Then we can use the probability measure to define  a conformal metric
\[
 ds_\rho^2=\rho(z)|\dd z|^2,\qquad
 \mathrm dA_\rho=\rho(z)\dA(z).
\]

\begin{proposition}\label{prop:curvature}
Put $r=|z|$ and $L=\log(1/r)$. As $z\to0$, uniformly in its argument,
\begin{align}
 \rho(z)&=\frac{2}{3\pi}\left(L-\frac16\right)+O(rL),
       \label{eq:sharpdensity}\\
 K_\rho(z)&=\frac{3\pi}{4r^2(L-\tfrac16)^3}\bigl(1+O(r)\bigr)
       \sim\frac{3\pi}{4r^2L^3}.\label{eq:curvature}
\end{align}
The same statements hold in the local coordinates $1-z$ at $1$ and
$1/z$ at $\infty$. In particular, the Gaussian curvature is positive
near each puncture and tends to $+\infty$.
\end{proposition}
\begin{proof}
The substitution in the proof of Proposition~\ref{prop:moments}, followed
by integration in $\beta$, gives
\[
 \rho(z)=\frac2\pi\int_0^1 V K(d,c)\dd v,
 \quad V=v(1-v),\quad d=r^2V,\quad
 c=1+(|1-z|^2-1)v,
\]
with the kernel $K$ from \eqref{eq:oned}. Splitting
$\log(c/d)=2L+\log c-\log V$ gives the local representation
$\rho=A L+B$ with $A,B$ real analytic across $0$.
Indeed, $c-d=|1-vz|^2$ stays uniformly away from zero near $0$,
and $V|\log V|$ is integrable. The integrals
\[ \int_0^1V\dd v=1/6, \quad and \quad \int_0^1V\log V\dd v=-5/18\] give
$A(0)=2/(3\pi)$ and $B(0)=-1/(9\pi)$, proving
\eqref{eq:sharpdensity} with control of its derivatives.

Let $b=B/A$, so $b(0)=-1/6$. For the Euclidean Laplacian $\Delta$,
\[
 K_\rho=-\frac{\Delta\log\rho}{2\rho},\qquad
 \Delta\log\rho=\Delta\log A+\frac{\Delta b}{L+b}
                   -\frac{|\nabla L+\nabla b|^2}{(L+b)^2}.
\]
Here $\Delta L=0$, $|\nabla L|=1/r$, and the derivatives of $A,b$
remain bounded. Thus
$\Delta\log\rho=-(1+O(r))/[r^2(L-1/6)^2]$, which proves
\eqref{eq:curvature}. Reflection and inversion give the other punctures.
\end{proof}

For comparison, write the curvature $-1$ hyperbolic metric as
$ds_{\mathrm{hyp}}=\lambda(z)|\dd z|$. Its cusp asymptotics are
\[
 \lambda(z)\sim\frac1{rL},\qquad
 \lambda(z)^2\sim\frac1{r^2L^2},
\]
see \cite{ASVV}, rescaling the curvature $-4$ convention there.
This shows the area densities to compare are $\rho$ and $\lambda^2$ and the first thing to observe is that
their singularities have different orders. More precisely, normalized
hyperbolic area has total mass one and satisfies
\[
 \rho\dA=F\dd\mu_{\mathrm{hyp}},\qquad
 \dd\mu_{\mathrm{hyp}}=\frac{\lambda^2}{2\pi}\dA,\qquad
 F=\frac{2\pi\rho}{\lambda^2}\sim\frac43r^2L^3\longrightarrow0.
\]
The weight $F$ is invariant under the anharmonic group.
Since $\int_0^\varepsilon\sqrt{\log(1/r)}\dd r<\infty$, each puncture is
at finite distance in $ds_\rho$. The hyperbolic distance integral
$\int_0^\varepsilon\dd r/[r\log(1/r)]$ diverges. Finally,
$r\partial_r\log\rho\to0$.
The positive curvature blow-up is integrable, with
\begin{equation}\label{eq:localcurvaturemass}
 \int_{0<|z|<\varepsilon}K_\rho\,\mathrm dA_\rho
 \sim\frac{\pi}{\log(1/\varepsilon)}\longrightarrow0.
\end{equation}

\section{Expectations and geometric quantities}
For any nonnegative measurable $h$, or any integrable $h$,
\begin{equation}\label{eq:fold}
 \E h(Z)=\int_{|z|<1}\bigl(h(z)+h(1/z)\bigr)\rho(z)\dA(z).
\end{equation}
If also $h(z)=h(1/z)=h(\ol z)$, this reduces to
\begin{equation}\label{eq:semidisc}
 \E h(Z)=4\int_{\substack{|z|<1\\\operatorname{Im}z>0}}
                   h(z)\rho(z)\dA(z).
\end{equation}
These are change-of-variable identities, with the Jacobian from
\eqref{eq:symmetry} included.

\subsection{Two exact logarithmic moments}
\begin{proposition}\label{prop:logmoments}
\begin{equation}\label{eq:logmoments}
 \E\log|Z|=0,\qquad \Var(\log|Z|)=\E(\log|Z|)^2=1.
\end{equation}
\end{proposition}
\begin{proof}
Let $q_{ij}=|P_i-P_j|/2$ be half the chordal distance between the corresponding
points on the unit sphere. The stereographic factors cancel in cross ratio,
so
\[
 \log|Z|=\log q_{13}+\log q_{24}-\log q_{12}-\log q_{34}.
\]
For independent uniform spherical points, $q_{ij}^2$ is uniform on $(0,1)$.
Thus each logarithm has mean $-1/2$ and variance $1/4$. Two different edges
with a common endpoint are conditionally independent given that endpoint,
and their conditional distributions do not depend on its position. They
therefore have covariance zero, as do disjoint edges. Taking the mean and
variance of the sum proves the claim. Joint independence of
all four logarithms does not hold.
\end{proof}

\subsection{Ideal tetrahedra}
The (unsigned) hyperbolic volume of the ideal tetrahedron determined by the
four points, in curvature $-1$, is $V=|D(Z)|$ \cite{Zagier}.
The expectation we obtain is a special case of
Theorem~1.1 of Kabluchko--Schange \cite{KS}: the convex hull of $n\geq4$
independent uniform points on the boundary sphere has expected volume
$\pi\bigl(n/2-\sum_{j=1}^{n-1}1/j\bigr)$.
We can give an evaluation from our density.
\begin{proposition}[Mean ideal-tetrahedron volume]\label{prop:volume}
\begin{equation}\label{eq:volume}
 \boxed{\displaystyle
 \int_\C |D(z)|\rho(z)\dA(z)=\E V=\frac{\pi}{6}.}
\end{equation}
\end{proposition}
\begin{proof}
Define the Lobachevsky function
\[
 \Lambda(\theta)=-\int_0^\theta\log(2\sin u)\dd u,
 \qquad 0\leq\theta\leq\pi.
\]
It is bounded and satisfies
$\Lambda(\pi-\theta)=-\Lambda(\theta)$.
If $\theta_1,\theta_2,\theta_3$ are the angles of the Euclidean triangle
with vertices $A,B,C$ in \eqref{eq:threepoints}, the usual angle formula
for ideal volume gives
$|D(Z)|=\sum_{j=1}^3\Lambda(\theta_j)$ \cite{Zagier}.
These three angles have the same distribution, since $A,B,C$ are
exchangeable. In particular, with $\Theta=|\arg Z|$ and
$m=\E\Lambda(\Theta)$, we have $\E V=3m$.

Put $g(\theta)=\int_0^\infty\rho(re^{i\theta})r\dd r$.
Conjugation symmetry and the antisymmetry of $\Lambda$ give
\[
 m=\int_0^\pi\Lambda(\theta)
       \bigl(g(\theta)-g(\pi-\theta)\bigr)\dd\theta.
\]
The radial integral in \eqref{eq:simplex} has an especially simple odd
part. Set
\[
 s=\frac{\alpha\beta}{(1-\alpha)(1-\beta)},\qquad
 k(s)=K(1,s)=\frac{-(1+s)\log s-2(1-s)}{(1-s)^3},
 \quad 0<s<1.
\]
Indeed, the substitutions $u=\alpha/(1-\alpha)$ and
$v=\beta/(1-\beta)$ turn the simplex into $u,v>0$, $uv<1$;
putting $v=s/u$ gives
\[
 \int_\Delta f(s)\dd\alpha\dd\beta
 =\int_0^1 f(s)\left(\int_0^\infty
       \frac{u\dd u}{(1+u)^2(u+s)^2}\right)\dd s
 =\int_0^1 f(s)k(s)\dd s.
\]
The last equality follows also from \eqref{eq:kernel}.
Rescaling the radial variable therefore gives
\[
 g(\theta)=\frac2\pi\int_0^1(1-s)k(s)
                      R(\sqrt{s}\cos\theta)\dd s,
 \qquad
 R(q)=\int_0^\infty\frac{r\dd r}{(r^2-2qr+1)^2}.
\]
The elementary identity
\[
 R(q)-R(-q)=\frac{\pi q}{2(1-q^2)^{3/2}},\qquad |q|<1,
\]
now yields
\[
 g(\theta)-g(\pi-\theta)
 =\cos\theta\int_0^1
  \frac{(1-s)\sqrt{s}\,k(s)}{(1-s\cos^2\theta)^{3/2}}\dd s.
\]
Integrating by parts in $\theta$, using
\[
 \frac{\mathrm d}{\mathrm d\theta}
 \frac{\sin\theta}{\sqrt{1-s\cos^2\theta}}
 =\frac{(1-s)\cos\theta}{(1-s\cos^2\theta)^{3/2}},
\]
reduces the expectation to
\begin{equation}\label{eq:volumeangular}
 m=\int_0^1 k(s)F(s)\dd s,\qquad
 F(s)=-a\log s-\Lambda(2a),\quad a=\arcsin\sqrt{s}.
\end{equation}
Here the inner integral is
\[
 \sqrt{s}\int_0^\pi
 \frac{\sin\theta\log(2\sin\theta)}
      {\sqrt{1-s\cos^2\theta}}\dd\theta
 =-a\log s-\Lambda(2a):
\]
use $t=\cos\theta$, then $\sqrt{s}\,t=\sin v$, and
$\sin^2a-\sin^2v=\sin(a+v)\sin(a-v)$.
The initial interchanges are justified by boundedness of $\Lambda$ and
normalization of $\rho$. After the integration by parts, use
$\sin\theta/\sqrt{1-s\cos^2\theta}\leq1$ and absolute integrability of
both $\log(2\sin\theta)$ and $\sqrt{s}\,k(s)$ on their respective intervals.

To evaluate \eqref{eq:volumeangular}, introduce
\[
 B(s)=\frac{1-s+s\log s}{(1-s)^2},\qquad
 C(s)=\frac{\log s}{1-s}.
\]
Then $k=-B'$, $C'=B/s$, and
\[
 F'(s)=-\frac{a}{s}
        +\frac{\log(16(1-s))}{2\sqrt{s(1-s)}}.
\]
Two integrations by parts give
\[
 m=\frac\pi2+\frac12\int_0^1
 \frac{\log s/(1-s)+B(s)\log(16(1-s))}
      {\sqrt{s(1-s)}}\dd s.
\]
The boundary terms use $B(0)=1$, $B(1)=1/2$,
$F(0)=F(1)=0$, and
$aC(s)\to0$ at $0$, $aC(s)\to-\pi/2$ at $1$.
More precisely, $F(s)=O(\sqrt{s})$ at $0$ and
$F(s)=O(\sqrt{1-s}\,|\log(1-s)|)$ at $1$;
the integrals are absolutely convergent.

Put $s=t^2/(1+t^2)$ and write
\[
 \ell(t)=\log\frac{t^2}{1+t^2},\qquad
 h(t)=\log\frac{16}{1+t^2}.
\]
We obtain
\[
 m=\frac\pi2+\int_0^\infty
       \bigl[\ell(t)+(1+t^2\ell(t))h(t)\bigr]\dd t.
\]
Finally, use the elementary logarithmic integrals
\[
 \int_0^\infty\ell(t)\dd t=-\pi,\qquad
 \int_0^\infty\frac{\log(1+t^2)}{t^2+v^2}\dd t
       =\frac\pi v\log(1+v)\quad(v>0),
\]
together with
\[
 1+t^2\ell(t)=\int_0^1\frac{u}{t^2+u}\dd u.
\]
The logarithmic identities follow by integration by parts or
differentiation in a positive scale parameter. Fubini is legitimate:
separating $h(t)=\log16-\log(1+t^2)$ gives two finite nonnegative
integrals. Consequently,
\[
 \int_0^\infty(1+t^2\ell(t))h(t)\dd t
 =\pi\int_0^1\sqrt{u}\log\frac4{1+\sqrt{u}}\dd u
 =2\pi\int_0^1v^2\log\frac4{1+v}\dd v.
\]
Thus the required density integral has the elementary evaluation
\begin{align}\label{eq:volumeoneD}
 \int_\C |D(z)|\rho(z)\dA(z)
 &=-\frac{3\pi}{2}+6\pi\int_0^1v^2\log\frac4{1+v}\dd v\nonumber\\
 &=-\frac{3\pi}{2}+6\pi\cdot\frac5{18}
 =\frac\pi6.
\end{align}
For the last step, integration by parts gives
$\int_0^1v^2\log(1+v)\dd v=\tfrac23\log2-\tfrac5{18}$.
\end{proof}

\subsection{Distance between two random geodesics}
Let $\ell_{14}$ and $\ell_{23}$ be the complete geodesics with endpoint pairs
$(Z_1,Z_4)$ and $(Z_2,Z_3)$. Their distance is
\begin{equation}\label{eq:distance}
 H(Z)=\arccosh\frac{1+|Z|}{|1-Z|}.
\end{equation}
To verify the pairing and formula, normalize the ordered quadruple to
$(0,1,z,\infty)$. The first geodesic is the vertical axis and the second
has endpoints $1,z$ in the upper half-space model. In general the minimum
distance $h$ from the vertical axis to a geodesic with finite endpoints
$u,v$ satisfies
\[
 \cosh^2h=\min_{-1<t<1}
 \frac{2(|u|^2+|v|^2)+2(|u|^2-|v|^2)t}
 {|u-v|^2(1-t^2)}
 =\frac{(|u|+|v|)^2}{|u-v|^2}.
\]
The first equality parametrizes the Euclidean semicircle defining the
geodesic; minimizing the quotient gives the second. Taking $u=1,v=z$
proves \eqref{eq:distance}. Other endpoint pairings have the same
distribution by exchangeability, though as functions of the
cross-ratio they differ.

$H$ satisfies the invariance needed in \eqref{eq:semidisc}.
It is integrable and its only unbounded behaviour is logarithmic near $z=1$,
where the density also has only a logarithmic singularity.
We now evaluate the expectation exactly. The useful observation is that
$H$ is harmonic off the negative real axis and its singularity at $1$.
Consequently its expectation can be recovered from a one-dimensional
trace of the logarithmic potential avoiding the closed dilogarithmic density.

\begin{lemma}[Logarithmic potential]\label{lem:potential}
Let $U(w)=\E\log|Z-w|$. For $w\in\C$, put
\[
 Q_w=\alpha\beta|1-w|^2+\alpha\gamma|w|^2+\beta\gamma.
\]
Then
\begin{equation}\label{eq:potential}
 U(w)=1+\int_\Delta\log Q_w\dd\alpha\dd\beta,
 \qquad U(0)=U(1)=0.
\end{equation}
In particular, for $t>0$,
\begin{equation}\label{eq:negativepotential}
 \boxed{\displaystyle
 U(-t)=-\frac13+
 \frac{(1+t)^3\log(1+t)-t^3\log t}{3t(1+t)}.}
\end{equation}
\end{lemma}
\begin{proof}
For $(\alpha,\beta)\in\Delta$, set
\[
 m=\frac{\beta}{1-\alpha},\qquad
 s^2=\frac{\beta\gamma}{\alpha(1-\alpha)^2}.
\]
Then $Q_w=\alpha(1-\alpha)(|w-m|^2+s^2)$, so
\eqref{eq:simplex} is a probability mixture:
\begin{equation}\label{eq:mixture}
 \rho(z)=2\int_\Delta
 \frac{s^2}{\pi(|z-m|^2+s^2)^2}\dd\alpha\dd\beta.
\end{equation}
The factor $2$ normalizes the area of $\Delta$.
The logarithmic potential of the individual density in this mixture is
\begin{equation}\label{eq:componentpotential}
 \int_\C\log|z-w|\,
 \frac{s^2\dA(z)}{\pi(|z-m|^2+s^2)^2}
 =\frac12\log(|w-m|^2+s^2).
\end{equation}
For example, the circle mean of $\log|re^{i\theta}-(w-m)|$ is
$\log\max\{r,|w-m|\}$; integrating this against
$2s^2r\dd r/(r^2+s^2)^2$ gives \eqref{eq:componentpotential}.
Taking the mixture and using
\[
 \int_\Delta\log[\alpha(1-\alpha)]\dd\alpha\dd\beta=-1
\]
proves \eqref{eq:potential}. These integrals are absolutely convergent:
the density has at most logarithmic singularities at $0,1$ and the tail
bound established in Proposition~\ref{prop:moments}.
For $w=0$, $Q_0=\beta(1-\beta)$ and its logarithm integrates to $-1$;
the same argument with $Q_1=\gamma(1-\gamma)$ gives $U(1)=0$.

For $w=-t$, introduce
\[
 q=\beta-t\alpha,\qquad r=\beta+t^2\alpha.
\]
Their Jacobian has absolute value $t(1+t)$, and
$Q_{-t}=r-q^2$. The image of $\Delta$ is the triangle with vertices
$(0,0),(-t,t^2),(1,1)$. Its upper boundary is
$r=(1-t)q+t$; the lower boundary is $r=-tq$ for $-t<q<0$
and $r=q$ for $0<q<1$. Integrate first in $r$, using
$F(v)=v(\log v-1)$ with $F(0)=0$. Define
\[
 G(L)=\int_0^L F\bigl(u(L-u)\bigr)\dd u
     =L^3\left(\frac{\log L}{3}-\frac49\right).
\]
The upper boundary contributes $G(1+t)$, and the two lower pieces
contribute $G(t)$ and $G(1)$. Hence
\[
 \int_\Delta\log Q_{-t}\dd\alpha\dd\beta
 =\frac{G(1+t)-G(t)-G(1)}{t(1+t)}.
\]
Substitution in \eqref{eq:potential} proves
\eqref{eq:negativepotential}.
\end{proof}

\begin{lemma}[Reduction to the negative real axis]\label{lem:Hpotential}
With the probability measure
\[
 \dd\nu(t)=\frac{\dd t}{\pi\sqrt t(1+t)},\qquad t>0,
\]
one has
\begin{equation}\label{eq:Hpotential}
 \E H(Z)=\int_0^\infty U(-t)\dd\nu(t).
\end{equation}
\end{lemma}
\begin{proof}
Choose the square root with positive real part on
$\C\setminus(-\infty,0]$. Directly from \eqref{eq:distance},
\begin{equation}\label{eq:Hharmonic}
 H(z)=\log\left|\frac{1+\sqrt z}{1-\sqrt z}\right|
     =2\log|1+\sqrt z|-\log|1-z|.
\end{equation}
Indeed,  the first logarithm is nonnegative and its hyperbolic cosine is
$(1+|z|)/|1-z|$. The elementary logarithmic-potential identity
\begin{equation}\label{eq:halfCauchy}
 2\log|1+\sqrt z|
 =\int_0^\infty\log|z+t|\dd\nu(t)
\end{equation}
follows by putting $t=u^2$ and differentiating in positive real $z$:
the derivative on both sides is $1/[\sqrt z(1+\sqrt z)]$.
The constant is fixed by the limit $z\downarrow0$, for which
$\int\log t\dd\nu(t)=0$ by $t\mapsto1/t$.
The corresponding holomorphic identity, with principal logarithms,
extends to the slit plane. Then taking real parts gives
\eqref{eq:halfCauchy}. Boundary values on the slit follow by continuity
and are immaterial for the expectation.

For completeness, Fubini's theorem applies absolutely here.
The positive part obeys
\[
 \log^+|Z+t|\leq\log(1+|Z|)+\log(1+t),
\]
whose expectation and $\nu$-integral are finite. For the negative part,
write $\log^-r=\max\{-\log r,0\}$. The density of $\nu$, extended by
zero to the negative half-line, belongs to $L^{3/2}(\R)$, whereas
$u\mapsto\log^-|u|$ belongs to $L^3(\R)$. Since
$\log^-|Z+t|\leq\log^-|\operatorname{Re}Z+t|$, H\"older's inequality
gives a bound on its $\nu$-integral uniform in $Z$.
Taking expectations in \eqref{eq:Hharmonic}--\eqref{eq:halfCauchy}
therefore yields
$\E H(Z)=\int U(-t)\dd\nu(t)-U(1)$.
Lemma~\ref{lem:potential} gives $U(1)=0$.
\end{proof}

\begin{theorem}[Mean distance between random geodesics]\label{thm:meandistance}
For two geodesics whose four endpoints are independent and uniform on the
boundary sphere, in curvature $-1$,
\begin{equation}\label{eq:meandistance}
 \boxed{\displaystyle
 4\int_{\substack{|z|<1\\\operatorname{Im}z>0}}
 H(z)\rho(z)\dA(z)=\E H(Z)=\frac43.}
\end{equation}
\end{theorem}
\begin{proof}
Insert \eqref{eq:negativepotential} in \eqref{eq:Hpotential} and put
$t=u^2$. This gives
\begin{align*}
 \E H(Z)=-\frac13+\frac{2}{3\pi}\int_0^\infty
 \biggl[&\log(1+u^{-2})+\frac{\log(1+u^2)}{u^2}\\
        &+\frac{2\log u^2}{1+u^2}
          -\frac{\log u^2}{(1+u^2)^2}\biggr]\dd u.
\end{align*}
Every term is individually integrable. The first two
integrals both equal $\pi$: substitution $u\mapsto1/u$ equates them,
and integration by parts evaluates the second. Inversion also gives
\[
 \int_0^\infty\frac{\log u^2}{1+u^2}\dd u=0.
\]
Finally, let $J=\int_0^\infty\log u^2\dd u/(1+u^2)^2$.
Inversion followed by addition and integration by parts yields
\[
 2J=\int_0^\infty\frac{(1-u^2)\log u^2}{(1+u^2)^2}\dd u
 =\int_0^\infty\log u^2\,\mathrm d\!\left(\frac{u}{1+u^2}\right)
 =-2\int_0^\infty\frac{\dd u}{1+u^2}=-\pi.
\]
All boundary terms vanish. Consequently
\[
 \E H(Z)=-\frac13+\frac{2}{3\pi}
 \left(\pi+\pi+0+\frac\pi2\right)=\frac43.
\]
Equation \eqref{eq:semidisc} supplies the semicircle integral.
\end{proof}

\section{Geometric applications}\label{sec:applications}

\subsection{Dihedral angles and acute ideal tetrahedra}
The mixture \eqref{eq:mixture} has a useful independent-variable form.
Writing $a=\alpha$ and $m=\beta/(1-\alpha)$ gives
$\beta=(1-a)m$, $\gamma=(1-a)(1-m)$ and
$\dd\alpha\dd\beta=(1-a)\dd a\dd m$. Consequently
\begin{equation}\label{eq:independentmixture}
 Z\overset d=M+\sqrt{\frac{M(1-M)}{A}}\,W,
\end{equation}
where $M$ is uniform on $(0,1)$, $A$ has density $2(1-a)$ on $(0,1)$,
and $W$ has the spherical density \eqref{eq:spherical}; the three variables
are independent.

Let $Z^+=\Re Z+i|\Im Z|$. The angles $\theta_1,\theta_2,\theta_3$ of
the Euclidean triangle with vertices $0,1,Z^+$, in that order, are the
three distinct interior dihedral angles of the corresponding ideal
tetrahedron. Opposite edges have equal angles and
$\theta_1+\theta_2+\theta_3=\pi$. The angles are exchangeable: in the
normalization \eqref{eq:threepoints}, they are the angles of a triangle
whose three vertices are independent with the same spherical law.

\begin{proposition}[Joint dihedral-angle density]\label{prop:angles}
For $u>0$, $v>0$, $w=\pi-u-v>0$, the joint density of
$(\theta_1,\theta_2)$ with respect to $\dd u\dd v$ is
\begin{equation}\label{eq:angledensity}
 p(u,v)=2\rho\!\left(\frac{\sin v}{\sin w}e^{iu}\right)
                   \frac{\sin u\sin v}{\sin^3w}.
\end{equation}
\end{proposition}
\begin{proof}
Conjugation symmetry gives $Z^+$ density $2\rho$ on the upper half-plane.
The sine rule gives $z=r e^{iu}$ with $r=\sin v/\sin w$.
Since $\partial_v r=\sin u/\sin^2w$, the planar Jacobian is
$r\partial_vr=\sin u\sin v/\sin^3w$. This parametrizes the upper
half-plane bijectively and proves both the formula and its normalization.
\end{proof}

\begin{proposition}[Acute ideal tetrahedra]\label{prop:acute}
For four independent uniform spherical vertices,
\begin{equation}\label{eq:acute}
 \boxed{\Prob(\hbox{all dihedral angles are acute})
       =\frac{3\pi^2}{8}-\frac72
       =0.2011\ldots .}
\end{equation}
\end{proposition}
\begin{proof}
Put $q=\Prob(\theta_1>\pi/2)=\Prob(\Re Z<0)$.
Integrating the spherical density over the  imaginary $i\mathbb{R}$ gives
\[
 f_{\Re W}(x)=\frac1{2(1+x^2)^{3/2}},\qquad
 \Prob(\Re W<t)=\frac12\left(1+\frac{t}{\sqrt{1+t^2}}\right).
\]
Conditioning in \eqref{eq:independentmixture} therefore yields
\begin{equation}\label{eq:obtuseintegral}
 q=\int_0^1\!\int_0^1(1-a)
       \left(1-\sqrt{\frac{am}{1-m+am}}\right)\dd m\dd a.
\end{equation}
These integrands are nonnegative and bounded by an integrable function.
For $0<a<1$, the substitution $(1-a)m=\sin^2t$ gives
\[
 \int_0^1\sqrt{\frac{am}{1-(1-a)m}}\dd m
 =\frac{\sqrt a}{(1-a)^{3/2}}
        \left(\arccos\sqrt a-\sqrt{a(1-a)}\right).
\]
Substitution into \eqref{eq:obtuseintegral}, then $a=\cos^2t$, gives
\[
 q=1-\int_0^1\sqrt{\frac a{1-a}}\arccos\sqrt a\dd a
   =1-2\int_0^{\pi/2}t\cos^2t\dd t
   =\frac32-\frac{\pi^2}{8}.
\]
At most one of the three distinct angles is obtuse, and right angles
have probability zero by absolute continuity. Exchangeability now gives
$\Prob(\hbox{all acute})=1-3q$, proving \eqref{eq:acute}.
\end{proof}

\subsection{Branched covers, lattices and punctured tori}\label{sec:modular}
The double cover of the sphere branched over four distinct points has
genus one: Riemann--Hurwitz gives $2g-2=2(-2)+4=0$.
After normalizing the branch points, it has the Legendre equation
\begin{equation}\label{eq:legendrecurve}
 E_z:\qquad v^2=u(u-1)(u-z).
\end{equation}
Choosing a branch point as origin makes this an elliptic curve; the
underlying conformal torus, and its once-punctured conformal class, are
independent of that choice. Forgetting the ordering of the branch points
identifies the six anharmonic values of $z$. Equivalently, their common
elliptic invariant is
\[
 j(z)=256\frac{(1-z+z^2)^3}{z^2(1-z)^2}.
\]

Let $\lambda$ be the modular lambda function and write
\[
 \mathcal F=\{\tau=x+iy:\ y>0,\ |x|\leq\tfrac12,\ |\tau|\geq1\}.
\]
Boundary identifications and the two elliptic points have area zero and
won't affect the probability. Let $\dd x\dd y$ be planar
area in  $\tau$.

\begin{proposition}[The modular probability law]\label{prop:modularlaw}
The conformal class of \eqref{eq:legendrecurve}, for $Z$ with density
$\rho$, has the following normalized density on $\mathcal F$:
\begin{equation}\label{eq:modulardensity}
 \boxed{\dd\nu(\tau)=\Phi(\tau)\dd x\dd y,\qquad
 \Phi(\tau)=6\rho(\lambda(\tau))|\lambda'(\tau)|^2.}
\end{equation}
Relative to normalized hyperbolic area,
\begin{equation}\label{eq:modularweight}
 \dd\mu_{\mathrm{mod}}=\frac3{\pi y^2}\dd x\dd y,\qquad
 \frac{\dd\nu}{\dd\mu_{\mathrm{mod}}}
       =2\pi y^2\rho(\lambda(\tau))|\lambda'(\tau)|^2.
\end{equation}
\end{proposition}
\begin{proof}
The map $\lambda$ identifies $\mathbb H/\Gamma(2)$ with
$\C\setminus\{0,1\}$; the subgroup $\Gamma(2)$ has index six in
$\mathrm{PSL}(2,\mathbb Z)$. Under the six cosets, $\lambda$ changes
by precisely the six anharmonic transformations \cite{DLMF}.
The transformation law \eqref{eq:symmetry} and the chain rule therefore
make $\rho(\lambda(\tau))|\lambda'(\tau)|^2\dd x\dd y$ invariant
under the full modular group. Its mass on a $\Gamma(2)$ fundamental
domain is one by change of variables. Such a domain consists, up to
null boundaries, of six translates of $\mathcal F$, so each has mass
$1/6$. This proves \eqref{eq:modulardensity}. Since the hyperbolic area
of $\mathcal F$ is $\pi/3$, division gives \eqref{eq:modularweight}.
\end{proof}

This  branched-cover construction for random tori has modulus obtained by pairing opposite sides of an ideal quadrilateral different to the real cross ratio \cite{Martin}. In
particular,  explicit elliptic-integral inversion of $\lambda$ does
not solve the accessory-parameter  for that side-pairing modulus.

For a torus represented by $\tau\in\mathcal F$, let $K_\square$ be the
least quasiconformal distortion of a map to the square torus, and let
$K_{\mathrm{rect}}$ be the least distortion of a map to any rectangular
torus. With the curvature $-1$ hyperbolic distance, the standard
extremal-distortion formulas \cite{MartinLattices} are
\begin{align}
 R_\square:=\log K_\square
 &=\arccosh\frac{x^2+y^2+1}{2y},\label{eq:squaredistortion}\\
 \log K_{\mathrm{rect}}&=\arsinh\frac{|x|}{y},\qquad
 K_{\mathrm{rect}}=\frac{\sqrt{x^2+y^2}+|x|}{y}
          \leq\sqrt3.\label{eq:rectdistortion}
\end{align}
Thus the conventional Teichm\"uller distance $\tfrac12\log K$ is
$R_\square/2$. We retain $\log K$ in all numerical comparisons.

For completeness, the square-torus representative $i$ is nearest among
its modular translates when $\tau\in\mathcal F$. Indeed, for a matrix
$\left(\begin{smallmatrix}a&b\\c&d\end{smallmatrix}\right)$ of determinant
one, put $n=c^2+d^2$, $m=a^2+b^2$ and $k=ac+bd$. The difference of the
two hyperbolic cosines, multiplied by $2y$, is
\[
 (n-1)|\tau|^2+(m-1)-2kx
 \geq n+m-2-|k|\geq0.
\]
Here $n,m\geq1$ and $n+m\geq|k|+2$: for $k\ne0$, use
$n+m-2|k|\geq1$ and integrality; for $k=0$ it is immediate.
The closest rectangular representative lies on the imaginary axis.
For $0\leq x\leq1/2$, the other two edges of the Farey triangle
$(0,1,\infty)$ have distances whose hyperbolic sines are
$(1-x)/y$ and $(|\tau|^2-x)/y$, both at least $x/y$.
Any path to another translate of the imaginary axis crosses an edge
of this triangle. Reflection handles $x<0$ and proves
\eqref{eq:rectdistortion}; its upper bound follows from
$|x|/y\leq1/\sqrt3$.

\begin{proposition}[Cusp tail and distortion moments]\label{prop:modulartail}
Under the law $\nu$, as $Y\to\infty$,
\begin{align}
 \Prob_\nu(\Im\tau>Y)
 &=512\left(\pi Y-\log16+\frac13\right)e^{-2\pi Y}
       +O(Ye^{-3\pi Y})\label{eq:modulartail}\\
 &\sim512\pi Y e^{-2\pi Y}.\nonumber
\end{align}
Moreover
\begin{equation}\label{eq:distortionmoments}
 \E_\nu K_\square^p<\infty\quad\hbox{for every }p>0.
\end{equation}
For normalized hyperbolic area and $p>0$, in contrast,
\begin{equation}\label{eq:uniformtail}
 \Prob_{\mu_{\mathrm{mod}}}(\Im\tau>Y)=\frac3{\pi Y}\quad(Y\geq1),
 \qquad
 \E_{\mu_{\mathrm{mod}}}K_\square^p<\infty
       \Longleftrightarrow 0<p<1.
\end{equation}
\end{proposition}
\begin{proof}
Writing $q=e^{\pi i\tau}$, the standard expansion is
$\lambda(\tau)=16q+O(q^2)$ and
$\lambda'(\tau)=16\pi i q+O(q^2)$ \cite[Eq.~23.17.4]{DLMF}.
These estimates are uniform for $|x|\leq1/2$. Substitution into
\eqref{eq:sharpdensity} and \eqref{eq:modulardensity} gives
\[
 \Phi(x+iy)=1024\pi\left(\pi y-\log16-\frac16\right)e^{-2\pi y}
                +O(ye^{-3\pi y}).
\]
For $Y\geq1$ the cusp is the rectangle $|x|\leq1/2$, $y>Y$.
Integrating first in $x$ and then in $y$ proves
\eqref{eq:modulartail}. Also \eqref{eq:squaredistortion} gives
\[
 K_\square=y+\frac{x^2}{y}+O(y^{-3})\qquad(y\to\infty)
\]
uniformly in $x$. Polynomial moments are therefore integrable against
the exponential tail, while the remaining part of $\mathcal F$ is
compact and $K_\square$ is bounded there. For hyperbolic area the same
comparison reduces convergence to $\int_1^\infty y^{p-2}\dd y$,
and integration of $3/(\pi y^2)$ gives the stated tail.
\end{proof}

The corresponding distribution functions have explicit integration
regions. For $r\geq0$,
\begin{align}
 \Prob_\nu(K_\square\leq e^r)
 &=\int_{\substack{\tau\in\mathcal F\\
        x^2+(y-\cosh r)^2\leq\sinh^2r}}
                      \Phi(\tau)\dd x\dd y,\label{eq:squarecdf}\\
 \Prob_\nu(K_{\mathrm{rect}}\leq e^r)
 &=\int_{\substack{\tau\in\mathcal F\\ |x|\leq y\sinh r}}
                      \Phi(\tau)\dd x\dd y.\label{eq:rectcdf}
\end{align}
These formulas also apply to hyperbolic area after replacing $\Phi$ by
$3/(\pi y^2)$.  The branched-cover law
has smaller mean distortion to the square torus but larger mean
distortion to the rectangular locus. Thin lattices lie close to that
locus while remaining far from the square torus, so these orderings
are compatible.

\section{Histograms and numerical checks}
\begin{landscape}
\begin{figure}[p]
\centering
\begin{minipage}[c]{0.32\linewidth}
\includegraphics[width=\linewidth]{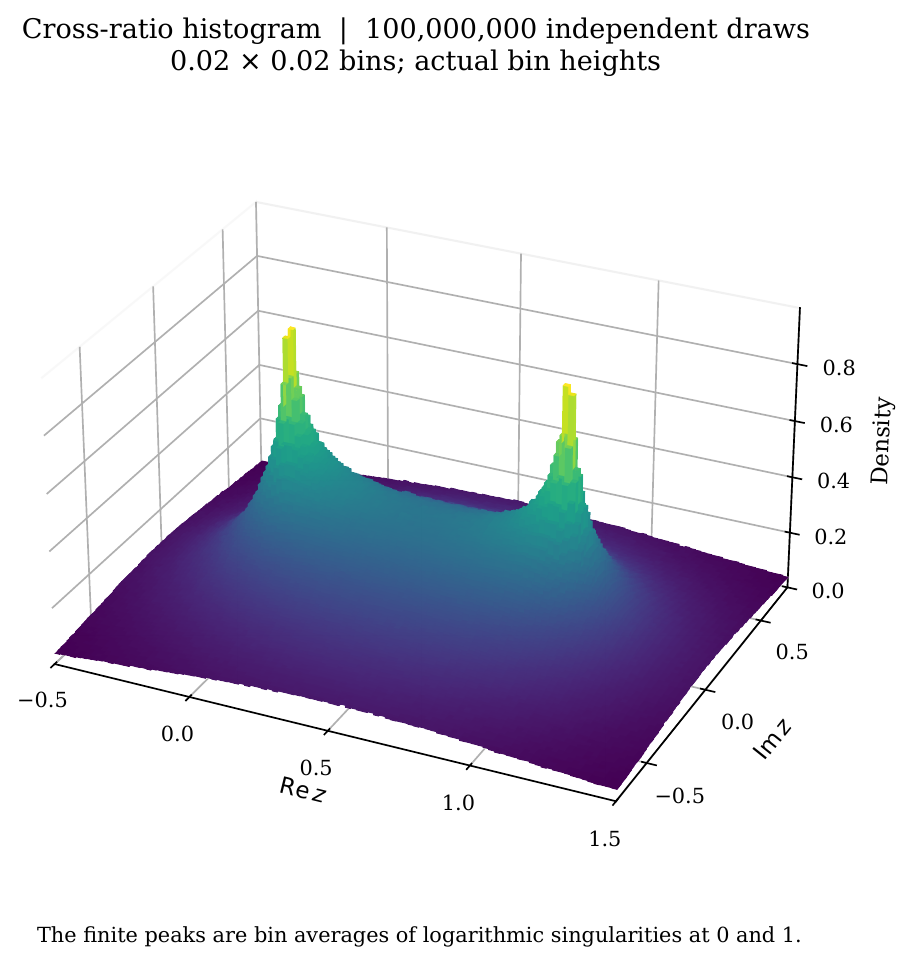}
\end{minipage}\hfill
\begin{minipage}[c]{0.66\linewidth}
\includegraphics[width=\linewidth]{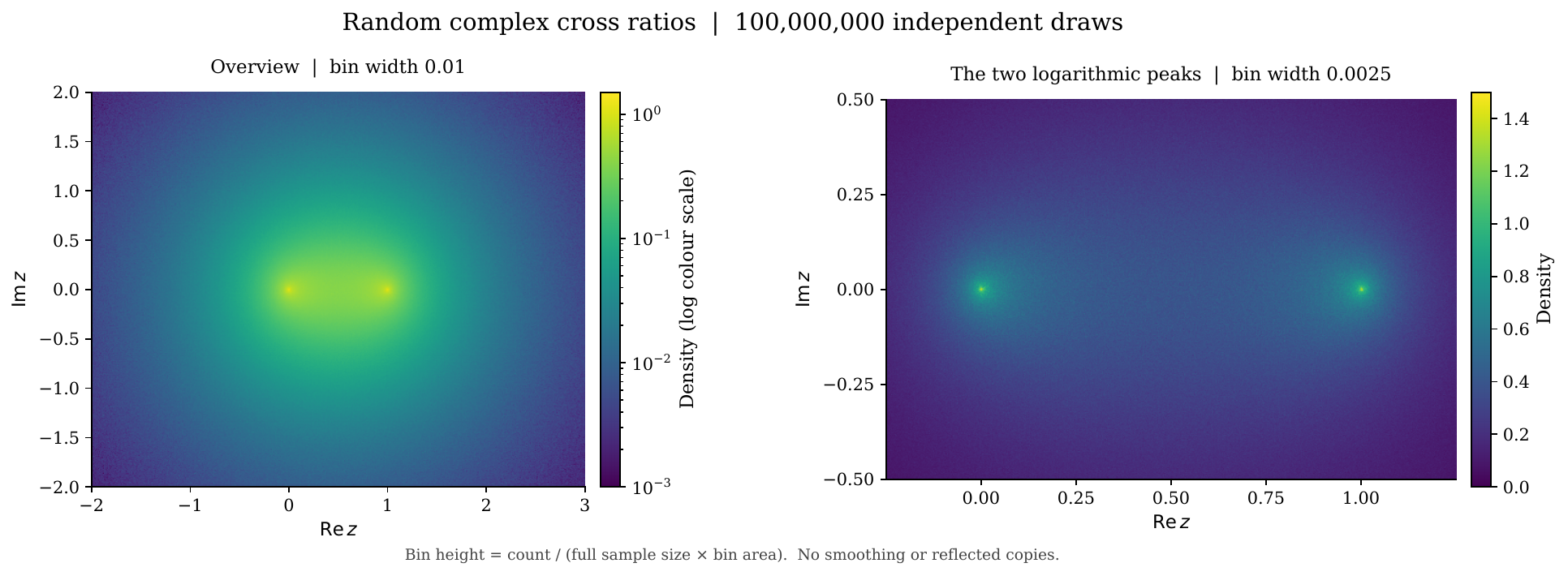}
\end{minipage}
\caption{Cross-ratio density from $100{,}000{,}000$ independent draws.
Left: three-dimensional histogram with square bins of side $0.02$.
Right: overview (left panel, bin side $0.01$, logarithmic colour scale)
and enlargement around $0$ and $1$ (right panel, bin side $0.0025$,
linear colour scale). All panels use the unconditional density
normalization \eqref{eq:histogram}, without smoothing or reflected copies.
The finite peaks are bin averages of the integrable logarithmic
singularities at $0$ and $1$.}
\label{fig:crossratio}
\end{figure}
\end{landscape}

\begin{figure}[!htbp]
\centering
\includegraphics[width=\textwidth]{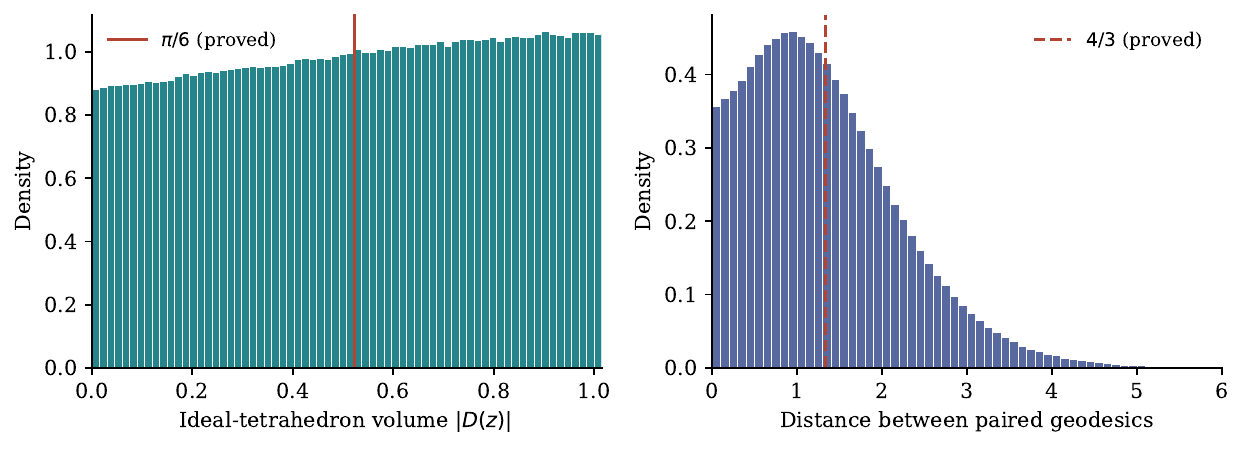}
\caption{The two-million-quadruple sample mapped to ideal-tetrahedron
volume and geodesic distance. The solid and dashed lines mark the proved volume and distance
means, respectively. The volume bins end at the regular
ideal tetrahedron's volume $D(e^{i\pi/3})=1.0149416064\ldots$.
Distances beyond the range were retained in all statistics and
in the histogram denominator.}
\label{fig:geometry}
\end{figure}

Figure~\ref{fig:crossratio} uses $N=10^10$ independent draws
from the exact mixture representation \eqref{eq:mixture}, with seed
$2026102501$. A separate sample of $2m$ independent quadruples,
generated using \eqref{eq:sampler} with seed $20261025$ yields  Figure~\ref{fig:geometry}.
Both experiments use NumPy's PCG64 generator. For a rectangular bin $B$
of area $|B|$, the plotted height is
\begin{equation}\label{eq:histogram}
 \widehat\rho_B=\frac{\#\{j:Z^{(j)}\in B\}}{N|B|}.
\end{equation}
Cropping the window does not change the denominator $N$.
No symmetry copies were added, and the bin heights were not smoothed.

The closed form was compared with the independent positive quadrature
\eqref{eq:oned} at fourteen points, including a point with imaginary part
$0.01$. The largest absolute discrepancy was below $10^{-9}$.
Reflection, conjugation and inversion identities were also checked
numerically; an exact symbolic calculation verified the differentiation
from \eqref{eq:rhooperator} to \eqref{eq:closed}.
 
\FloatBarrier

The ordinary sample mean of $Z$ is a poor diagnostic despite its finite
expectation, because its second absolute moment is infinite
(Proposition~\ref{prop:moments}). The bounded symmetry indicators and the
logarithmic moments provide better checks of the sampling law. 
 
  \begin{figure}[!htbp]
\centering
\includegraphics[width=\textwidth]{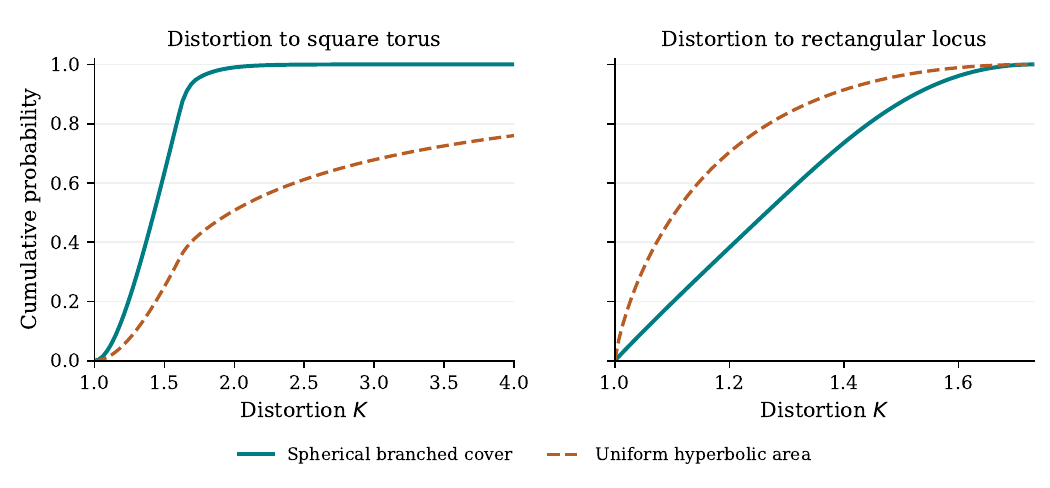}
\caption{Cumulative distributions of the optimal quasiconformal
distortion to the square torus (left) and to the rectangular locus
(right). Solid curves use the spherical branched-cover law
\eqref{eq:modulardensity}; dashed curves use normalized hyperbolic area.
The curves are numerical integrals over the regions
\eqref{eq:squarecdf}--\eqref{eq:rectcdf}. The rectangular distortion is
bounded by $\sqrt3$ for both laws.}
\label{fig:modular}
\end{figure}
\FloatBarrier

\end{document}